\documentclass[12pt]{article}
\usepackage{amssymb}
\newtheorem{theorem}{\sc Theorem}[section]
\newtheorem{definition}{\sc Definition}[section]
\newtheorem{lemma}{\sc Lemma}[section]

\newtheorem{example}{\sc Example}[section]
\newtheorem{remark}{\it Remark}[section]

\begin{document}
\baselineskip=24 pt
\title{\bf Integral manifolds of  differential equations with piecewise constant
argument of generalized type}
\author{ M. U. Akhmet\thanks{M.U. Akhmet is previously known as M. U. Akhmetov. }}
\date{{\small Department of Mathematics and Institute of Applied Mathematics, Middle East
Technical University, 06531 Ankara, Turkey}}
\maketitle

\vspace{0.5cm}


\begin{abstract}
In this  paper we introduce a general type of differential equations with piecewise constant argument (EPCAG), and consider the problem of backward continuation of solutions. We establish the existence of global integral manifolds of quasilinear EPCAG, consisting of solutions back continued  to $-\infty$, while the solutions  starting outside the invariant sets may not be back continued. The smoothness of the manifolds is investigated.  The existence of bounded and periodic solutions is considered. A new technique of investigation of equations with piecewise argument, based on an integral representation formula, is proposed. 
\end{abstract}
\maketitle
\section{Introduction and Preliminaries}
\subsection{Definitions and the description of the system}
Let $\mathbb Z, \mathbb N$ and $\mathbb R$ be the  sets of all integers, natural and real
numbers, respectively. Denote by $||\cdot||$ the Euclidean norm in
$\mathbb R^n$, $n \in \mathbb N.$

In this paper we are concerned with the quasilinear system
\begin{eqnarray}
&& y' = A(t)y + f(t,y(t),y(\beta(t))),
\label{1}
\end{eqnarray}
where $y \in \mathbb R^n, t \in \mathbb R, \beta(t) =\theta_i$ 
if $\theta_i \leq t <\theta_{i+1}, i= \mathbb Z,$ is an identification function, 
$\theta_i, i \in \mathbb Z,$ is a strictly ordered sequence of real numbers,
$|\theta_i| \rightarrow \infty $ as $ |i| \rightarrow \infty,$ and there
exists a number $\theta >0$ such that $\theta_{i+1}-\theta_{i} \leq \theta,\, i \in \mathbb Z.$  
The theory of differential equations with  piecewise constant argument
 (EPCA) of the type
\begin{eqnarray}
&& \frac{dx(t)}{dt} = f(t,x(t),x([ t])),
\label{e2}
\end {eqnarray}
where $[\cdot]$ signifies the greatest integer function,  was initiated in  \cite{cw1} and 
developed by many authors
 \cite{aw1,aa1, cw1, ky,p,p1, s2}, \cite{wl}-\cite{y1}.

The novel idea of this paper is that 
 system (\ref{1})  is a general case (EPCAG)  of equation (\ref{e2}). 
Indeed, if we take $\theta_i = i, i  \in \mathbb Z,$ then   (\ref{1}) takes  the form of (\ref{e2}).

The existing method of investigation  of EPCA, as   proposed 
by its founders, is based on the reduction of EPCA to  discrete equations. 
We propose another  approach to the problem. In fact, this approach consists
of the construction of the equivalent integral 
equation. Consequently, for every result of our paper we prove  a corresponding equivalence lemma.
Thus, while investigating EPCAG, we need not impose any conditions on  the reduced discrete equations, and, hence,  we require more easily verifiable conditions, similar to those for ordinary differential equations. 
It become less cumbersome  to solve the problems  of EPCAG theory (as well as of   EPCA theory). 
 
The theory of integral manifolds was founded by  H. Poincar\'e and A. M. Lyapunov \cite{poin,lya},
and it became a very powerful instrument for investigating various problems of the qualitative theory 
of differential equations. For example, one can talk about the exceptional role of  manifolds    
in the reduction of the  dimensions of equations. It is natural that the exploration of manifolds, of their
properties and neighborhoods, is one of the most interesting problems 
\cite{b1, b, carr, h, kelley, pal, pal1, pliss, pliss1, pugh}. 
One should  not be surprised that 
 manifolds are  one of the major subjects of investigation for specific types of differential  and difference equations \cite{akhmet,akhmet1, h1, lbs, sp, stokes}.  EPCA are no exception \cite{p1}.  Obviously, it is not possible to mention all the results 
pertaining to  integral sets  in this paper.

In what follows, we use the uniform norm $||T|| = \sup\{||Tx|| | ||x||=1\}$ for matrices.
The following assumptions will be needed throughout the paper.
\begin{itemize}
\item[C1)]  $A(t)$ is a continuous $n\times n$ matrix and $\sup_{\mathbb R} ||A(t)|| = \mu < \infty;$
\item[C2)] $f(t,x,z)$ is continuous in the first argument, $f(t,0,0) = 0, t \in \mathbb R,$ and f is Lipshitzian in the second and  the third arguments with a Lipshitz
constant  $l$ such that
$$||f(t,y_1,w_1) -  f(t,y_2,w_2)|| \leq l ( ||y_1-y_2|| + ||w_1-w_2||);$$
\item[C3)] the linear homogeneous system associated with (\ref{1}) 
\begin{eqnarray}
&& x' = A(t)x
\label{2}
\end{eqnarray}
has an exponential dichotomy on $\mathbb R$. That is, there exists a projection $P$ and positive constants
$K \,\mbox{and}\, \sigma$ such that 
\begin{eqnarray}
&&||X(t)PX^{-1}(s)|| \leq K \exp(- \sigma (t-s)), \, t \geq s,\nonumber\\
&&||X(t)(I-P)X^{-1}(s)|| \leq K \exp( \sigma (s-t)), \, t \leq s,
\label{3}
\end{eqnarray}
where $X(t)$ is a fundamental matrix of (\ref{2}).
\end{itemize}
\begin{definition}  A function $x(t)=x(t,\theta_i,x_0), x(\theta_i) = x_0,$ is a solution of equation (\ref{1})
on the interval $[\theta_i, \infty)$ for a fixed $i \in \mathbb Z,$  if
the following conditions are fulfilled:
\begin{enumerate}
\item[(i)] $x(t)$ is continuous on $[\theta_i, \infty);$
\item[(ii)] the derivative $x'(t)$ exists at each point $t \in [\theta_i, \infty)$ with the possible exception
of the points $\theta_i \in [t_0, \infty),$  where the right-sided derivative exists;
\item[(iii)] equation (\ref{1}) is satisfied by $x(t)$ on each interval $(\theta_i, \theta_{i+1}), i \in \mathbb Z,$ and it holds for the right derivative of 
$x(t)$ at the points $t= \theta_i, i \in \mathbb Z.$ 
\end{enumerate}
\label{defn8}
\end{definition}
\begin{remark} One can see that Definition \ref{defn8} is a slightly  changed version of a definition from
\cite{cw1},  adapted for our general case. Moreover, we have specified unbounded intervals 
$[\theta_i, \infty)$ since one can easily check that every solution $x(t, \theta_i,x_0), x_0 \in \mathbb R^n, i \in \mathbb Z,$ of (\ref{1}) can be continued to $\infty.$
\end{remark}
\begin{definition} A point  $(t_0,x_0), x_0 \in \mathbb R^n, \theta_i < t_0 \leq \theta_{i+1},$ is said to be backward continued to $t=\theta_i$
if there exists a solution $x(t,\theta_i,\bar x), \bar x \in \mathbb R,$  of the equation 
$x' = A(t)x + f(t,x(t),\bar x)),$ 
such that $x(t_0,\theta_i,\bar x) =  x_0.$ If the continuation is unique, $ \quad(t_0,x_0)$ is uniquely backward continued to  $t=\theta_i$  . 
If the point $(t_0,x_0)$ is backward continued to $t=\theta_i,$ and $x(t,\theta_i,\bar x)$ is its
continuation, then we let $x(t,t_0,x_0) = x(t,\theta_i,\bar x)$ to be  a solution of  (\ref{1}) with the initial data  $(t_0,x_0),$ and shall say that  $x(t,t_0,x_0)$ is back continued to 
 $t=\theta_i.$
\label{defn2}
\end{definition}
The following  example  shows that even for simple EPCAG the backward continuation of some  solutions can fail. 
\begin{example} 
Consider the following EPCA	
\begin{eqnarray}
&& x'= 2x - x^2([t]),
\label{ex1} 	
\end{eqnarray}
where $x \in \mathbb R,  t \in \mathbb R.$ 
Let us show that not all solutions of (\ref{ex1}) can be continued back.
Consider the interval $[0,1].$  Fix $z \in \mathbb R$ and let
 $x(t,1,z)$ be a solution of  (\ref{ex1}). Introduce an operator
$T:   \mathbb R \rightarrow \mathbb R$ such that
$$T x_0 = \exp(2)x_0+ \int_0^1\exp(2(1-s))x_0^2ds.$$ It can be easily seen that for the solution
  $x(t,1,z)$ to be back continued to $t=0$ the equation
$T x_0 = z$ must be solvable with respect to $x_0.$ But the last equation is equivalent to the equation  $[\exp(2)-1]x_0^2 +2\exp(2)x_0 - 2z = 0.$ Since the latter  can be solved not for all $z \in \mathbb R,$ the assertion is proved. 
Further we shall consider the uniqueness of the continuation.
Consider the interval $[0,1]$ again.
Fix numbers $x_0, x_1 \in \mathbb R$ such that 
$(x_0 + x_1)(1-\exp(2)) = 2\exp(2).$ 
Denote $x_0(t) = x(t,0,x_0)$ and $ x_1(t) = x(t,0,x_1),$ solutions of (\ref{ex1}).
For $t \in [0,1],$ they are solutions of  the equations 
$x'= 2x - x_0^2$ and $x'= 2x - x_1^2,$ respectively. 
Since 
$x_j(t) = \exp(2t)x_j + \int_0^t\exp(2(t-s))x_j^2ds, j=0,1,$ one can find that
$x_0(1) =x_1(1).$ That is, the solution $x(t,1,x_1(1))$ of (\ref{ex1}) cannot be continued to
$t=0$ uniquely. 
\end{example}
\begin{definition} A function $x(t)=x(t,t_0,x_0), x(t_0) = x_0,  \theta_i < t_0 < \theta_{i+1}, i \in \mathbb Z,$
 is a solution of (\ref{1})
on the interval $[\theta_i, \infty)$  if
the following conditions are fulfilled:
\begin{enumerate}
\item[(i)] the point  $(t_0,x_0)$ is backward continued to $t=\theta_i;$
\item[(ii)] the derivative $x'(t)$ exists at each point $t \in [\theta_i, \infty)$ with the possible exception
of the points $\,\theta_j \in [\theta_i, \infty),$   where the right-sided derivative exists;
\item[(iii)] equation (\ref{1}) is satisfied at each point $t \in [\theta_i, \infty)\backslash \{\theta_i\},$ and it holds for the  right derivative of $x(t)$ at the points $\theta_j \in [\theta_i, \infty).$

\end{enumerate}
\label{defn3}
\end{definition}
\begin{definition} The solution $x(t)=x(t,t_0,x_0)$ of (\ref{1}) is said to be backward continued to $t=a, a < t_0,$
if there exists a solution of (\ref{1}) $x(t,a,\bar x),  \bar x \in \mathbb R^n,$ in the sense 
of  Definition  \ref{defn3}   such that $x(t_0,a,\bar x) = x_0.$ If the continuation is unique, $x(t)$ is uniquely backward
continued to  $t=a$  .
\label{defn14}
\end{definition}
\begin{definition} The solution $x(t)=x(t,t_0,x_0)$ of (\ref{1}) is said to be backward continued to $-\infty$
if it  is  backward continued to every $a \in \mathbb R, a < t_0,$  and
$x(t)$ is uniquely backward continued to  $-\infty$  if the continuation is unique.
\label{defn6}
\end{definition}
\begin{remark} The backward continuation of the solutions of EPCA was considered in \cite{cw1} through the solvability of  certain difference equations. For our needs,  we have introduced less formal definitions since we shall consider  integral manifolds,  and  it is natural  to discuss the global backward continuation as well as its uniqueness on these manifolds. The definition of backward continuation for functional differential equations is considered in \cite{h1}.

\end{remark}
\begin{definition} A solution  of (\ref{1}), $x(t)=x(t,t_0,x_0),$  is said to be continued 
on  $\mathbb R$ if it is continued to $\infty$ and backward  continued to $-\infty.$
\end{definition}
\begin{definition} 
The set $\Sigma$  in the $(t,x)-$ space is said to be  an  integral set of system (\ref{1}) if any solution
$x(t) = x(t,t_0,x_0), x(t_0) = x_0,$ with $(t_0,x_0) \in \Sigma,$ has the property that $(t,x(t)) \in \Sigma, t \geq t_0,$
and the solution is backward continued to $-\infty$ so that $(t,x(t)) \in \Sigma, t < t_0.$ In other words, for every $(t_0,x_0) \in \Sigma$ the solution 
$x(t) = x(t,t_0,x_0), x(t_0) = x_0,$ is continuable on  $\mathbb R$ and 
$(t,x(t)) \in \Sigma, t \in \mathbb R.$
\label{defn13}
\end{definition}
We shall also use  the following definition, which is a version of a definition from \cite{p1}, adapted for our general case.
\begin{definition}   A function $x(t)$ is a solution of (\ref{1}) on $\mathbb R$ if:
\begin{enumerate}
\item[(i)] $x(t)$ is continuous on $\mathbb R;$
\item[(ii)] the derivative $x'(t)$ exists at each point $t \in \mathbb R$ with the possible exception
of the points $\theta_i , i \in \mathbb Z,$  where the right-sided derivatives exist;
\item[(iii)] equation (\ref{1})  is satisfied on each interval $(\theta_i, \theta_{i+1}), i \in \mathbb Z,$ and it holds for the  right derivative of $x(t)$ at the points $\theta_i , i \in \mathbb Z.$ 

\end{enumerate}
\label{defn7}
\end{definition}
\begin{remark} It is obvious that every solution $x(t)=x(t,t_0,x(t_0)), t_0 \in \mathbb R,$ is backward continued to $- \infty $  if   $x(t)$ is a solution by Definition  \ref{defn7}.
Conversely, if a solution $x(t) = x(t,t_0,x_0), t \geq t_0,$ is backward continued to $-\infty,$ then, denoting the continuation
by $x(t)$ for  $t < t_0,$ one can see that the function $x(t)$ is a  solution 
in the sense of Definition \ref{defn7}.
The backward continuation  we have introduced in this paper is convenient  
for keeping  the similarity with the  definition of integral sets for ordinary differential equations. 
Moreover, the uniqueness of backward continuations is essential for proving the invariance of the surfaces.
\end{remark}
\subsection{The uniqueness of the backward continuation}
Since the integral manifolds which  we are going to consider in the next sections are invariant,
it is useful to investigate the problem of the uniqueness of the  backward continuation of solutions of equation (\ref{1}).
We shall use the following assertion from \cite{hw}.
\begin{lemma} Assume that  condition $C1)$ is fulfilled. Then
\begin{eqnarray*}
&& ||X(t,s)|| \leq \exp(\mu|t-s|), t,s \in \mathbb R.
\label{w1}
\end{eqnarray*}
\label{wlemma}
\end{lemma}
\begin{lemma} Assume that condition $C1)$ is fulfilled. Then
\begin{eqnarray*}
&& ||X(t,s)|| \geq  \, \exp(-\mu|t-s|), t,s \in \mathbb R.
\label{w2}
\end{eqnarray*}
\label{wlemma1}
 \end{lemma}
{\it Proof.} The proof follows immediately from the equality $X(t,s) X(s,t) = I,$
where $I$ is  an $n \times n$ identity   matrix. 

The last two lemmas imply  the following,  simple but useful in what follows, inequalities
\begin{eqnarray}
&& ||X(t,s)|| \leq M, \nonumber\\
&& ||X(t,s)|| \geq m,
\label{w3}
\end{eqnarray}
which hold if $|t-s| \leq \theta,$ where $M =  \exp(\mu \theta), m = \exp(-\mu \theta).$

The following assertion can be easily proved.

\begin{theorem} 
Assume that conditions $C1)-C3)$ are fulfilled, 
\begin{eqnarray}
&& l M \theta[1+ M(1+l \theta) \exp(Ml\theta)] < m.
\label{bw}
\end{eqnarray}
Then  every solution of  (\ref{1}) has an unique continuation.
\end{theorem}
\subsection{Reduction to a system with a box-diagonal matrix of coefficients}
Using Gram-Schmidt orthogonalization of the columns of $X(t)$  \cite{cop,p1},  one can obtain that  by the transformation  $y = U(t)z,$ where $U(t)$ is a Lyapunov  matrix, (\ref{1}) can be reduced to the following system
\begin{eqnarray}\nonumber
\frac{du}{dt} &=&B_+(t)u + g_+(t,z(t),z(\beta(t))), \nonumber\\
\frac{dv}{dt} &=&B_-(t)v + g_-(t,z(t),z(\beta(t))), 
\label{3*}
\end{eqnarray}
where $$z=(u,v), u \in {R}^k , v \in {R} ^{n-k}, diag\{B_+(t), B_-(t)\} = U^{-1}(t)A(t) U(t),$$ 
$$(g_+(t,z(t),z(\beta(t))),g_-(t,z(t),z(\beta(t)))) =f(t,U(t)z(t)), U(\beta(t))z(\beta(t))).$$
One can check that the Lipshitz condition is valid 
$$||g_+(t,z_1,w_1) - g_+(t,z_2,w_2) || + ||g_-(t,z_1,w_1) - g_-(t,z_2,w_2)|| $$ $$ \leq L (||z_1-z_2|| +||w_1-w_2||) $$
for all $t \in \mathbb R, z_1,z_2 \in \mathbb R^k, w_1,w_2 \in  \mathbb R^{(n-k)},$ and 
$L = 2\sup_{\mathbb R}||U(t)|| l.$

The normed fundamental matrices $U(t,s), V(t,s)$ of the systems 
\begin{eqnarray}
\frac{du}{dt} &=&B_+(t)u, \nonumber\\
\frac{dv}{dt} &=&B_-(t)v, 
\label{4}
\end{eqnarray}
respectively, satisfy the following inequalities
\begin{eqnarray}
&&||U(t,s)|| \leq K \exp(- \sigma (t-s)), \, t \geq s,\nonumber\\
&&||V(t,s)|| \leq K \exp( \sigma (s-t)), \, t \leq s.
\label{5}
\end{eqnarray}

\section{Main Results}
\subsection{The existence of manifolds}
The following two lemmas are of major importance for our paper and they can be verified by differentiation.
\begin{lemma} Fix $ N \in \mathbb R, N >0, \alpha \in (0,\sigma)$ and assume that conditions $C1)-C3)$
and inequality (\ref{bw}) are valid.  
A function $z(t) = (u,v), ||z(t)|| \leq N \exp (-\alpha (t-t_0)),
t \geq t_0,$ is a solution of (\ref{3*}) on $\mathbb R$ if and only if it is a solution
on $\mathbb R$ of the following
system of integral equations
\begin{eqnarray}
&& u(t) = U(t,t_0)u(t_0) + \int^t_{t_0}U(t,s)g_+(s,z(s),z(\beta(s)))ds,\nonumber\\
&& v(t) = - \int^{\infty}_t V(t,s)g_- (s,z(s),z(\beta(s)))ds.
\label{6}
\end{eqnarray}
\label{l1}
\end{lemma}
\begin{lemma} Fix $ N \in \mathbb R, N >0, \alpha \in (0,\sigma),$ and assume that conditions $C1)-C3)$
and inequality (\ref{bw}) are valid.  
A function $z(t) = (u,v), ||z(t)|| \leq N \exp (\alpha (t-t_0)),
t \leq t_0,$ is a solution of (\ref{3}) on $\mathbb R$ if and only if it is a solution of the following
system of integral equations
\begin{eqnarray}
&& u(t) =  \int_{-\infty}^{t}U(t,s)g_+(s,z(s),z(\beta(s)))ds,\nonumber\\
&& v(t) = V(t,t_0)v(t_0) + \int_{t_0}^t V(t,s)g_- (s,z(s),z(\beta(s)))ds.
\label{7*}
\end{eqnarray}
\label{l2}
\end{lemma}

The proof of the next theorems is  very similar to that of the classic assertions about integral manifolds 
\cite{h, pal, pal1, pliss, pliss1}.
\begin{theorem} Suppose that conditions $C1)-C3)$ and inequality  (\ref{bw}) are satisfied. Then  for arbitrary $\epsilon >0, \alpha \in (0, \sigma)$ and 
a sufficiently small Lipshitz
constant $L,$ there 
exists a continuous function
$F(t,u)$  satisfying
\begin{eqnarray}
&& F(t,0)=0,\\
&& ||F(t,u_1) - F(t,u_2)|| \leq \frac{2K^2L(1+\exp(\sigma\theta))}{\sigma +\alpha} ||u_1-u_2||,
\end{eqnarray}
for all $t, u_1, u_2,$ such that $v_0=F(t_0,u_0)$ determines  a solution $z(t)$ of  (\ref{3*})
which is continued on $\mathbb R$ and 
\begin{eqnarray}
&& ||z(t)|| \leq (K+ \epsilon) ||u_0|| \exp(-\alpha(t-t_0)), t \geq t_0.
\label{7***}
\end{eqnarray}
\label{t1}
\end{theorem}
{\it Proof.} Let us consider  system (\ref{6}) and apply the method of successive approximations to it.
Denote $z_0=(0,0)^T, z_m = (u_m,v_m)^T, m \in \mathbb N,$ where for $ m \geq 0$
\begin{eqnarray}
&& u_{m+1}(t) = U(t,t_0)u(t_0) + \int^t_{t_0}U(t,s)g_+(s,z_m(s),z_m(\beta(s)))ds,\nonumber\\
&& v_{m+1}(t) = - \int^{\infty}_t V(t,s)g_- (s,z_m(s),z_m(\beta(s)))ds.
\label{8}
\end{eqnarray}
One can show by induction  that 
\begin{eqnarray}
&& ||z_m(t,c)|| \leq (K+ \epsilon) ||c|| \exp(-\alpha(t-t_0)), t \geq t_0,
\label{9}
\end{eqnarray}
provided that 
\begin{eqnarray}
&& K(K+\epsilon) \frac{2\sigma}{\sigma^2-\alpha^2}(1+ \exp(\sigma\theta)) L < \epsilon.
\label{10}
\end{eqnarray}
Similarly, one can establish  the following inequalities 
\begin{eqnarray*}
&& ||v_m(t,c_1) - v_m(t,c_2) || \leq \frac{2K^2L(1+\exp(\sigma\theta))}{\sigma +\alpha} ||c_1-c_2|| \exp(-\alpha(t-t_0)),\\
&& ||z_m(t,c_1) - z_m(t,c_2) || \leq 2K ||c_1-c_2|| \exp(-\alpha(t-t_0)),
\end{eqnarray*}
if  $4\sigma KL(1+\exp(\sigma\theta))< \sigma^2-\alpha^2.$

 And 
\begin{eqnarray}
&&  ||z_m(t,c) - z_{m-1}(t,c) || \leq \nonumber\\
&& K ||c|| \Big( \frac{2KL(1+\exp(\sigma\theta))}{ \sigma - \alpha} \Big)^{m-1} \exp(-\alpha(t-t_0)).
\end{eqnarray}
The last inequality  and the assumption
\begin{eqnarray}
&& L < \frac{\sigma - \alpha}{2K(1+\exp(\sigma\theta))}
\label{12}
\end{eqnarray}
 imply that the sequence $z_m$ converges uniformly
for all $c$ and $t\geq t_0.$ Define  the limit function $z(t,t_0,c)=(u(t,t_0,c), v(t,t_0,c)).$
It can be  easily seen that this function is a solution of (\ref{6}). By Lemma \ref{l1} $z(t,t_0,c)$ is a solution of  (\ref{3*}), 
too.  Taking $t=t_0$ in (\ref{6}) we  have that
\begin{eqnarray*}
&& u(t_0,t_0,c) = c,\nonumber\\
&&	v(t_0,t_0,c) = - \int^{\infty}_{t_0} V(t,s)g_- (s,z(s,t_0,c),z(\beta(s),t_0,c)))ds.
\end{eqnarray*}

Denote $F(t_0,c) = v(t_0,t_0,c).$ One can see that it satisfies all the conditions which should be verified.
The Theorem is proved.

Let us denote by $S^+$ the set of all points from the $(t,z)-$ space such that $v=F(t,u).$
\begin{theorem} The set $S^+$ is an integral surface.
\label{tS}
\end{theorem}
{\it Proof.}  Assume that $(t_0,u_0,v_0) \in S^+, z_0 = (u_0,v_0).$ We must show that if
$z(t) = z(t, t_0,z_0),$ then $(t^*,z(t^*)) \in S^+$ for  all $t^* \in \mathbb R.$
Indeed, if $t^* > t_0,$ then $||z(t)|| \leq  (K+ \epsilon) ||u_0||  \exp(-\alpha(t^*-t_0))\exp(-\alpha(t-t^*)), t^* > t_0.\,$
Lemma \ref{l1} implies that the point $(t^*,z(t^*))$ satisfies the equation $v=F(t,u).$
If  $t^* < t_0,$ then $||z(t)|| \leq  \bar K||u_0|| \exp(-\alpha(t^*-t_0))\exp(-\alpha(t-t^*)),$
where $\bar K = \max\{\max_{[t^*,t_0]}||z(t)||,  (K+ \epsilon) ||u_0||\}.$ Using Lemma  \ref{l1} again one can see
that $(t^*,z(t^*)) \in S^+.$
The Theorem is proved.
\begin{theorem} For every  fixed   $(t_0, u_0)$  system (\ref{6}) admits only one solution bounded
 on $[t_0, \infty).$
\label{t2}
\end{theorem}
{\it Proof.} \rm If  $z_1, z_2$ are two bounded solutions of  (\ref{6}),  then by straightforward evaluation it can be shown that
$$ \sup_{[t_0,\infty)} ||  z_1 -  z_2|| \leq  \frac{2KL(1+\exp(\sigma))}{\sigma - \alpha} \sup_{[t_0,\infty)} ||
  z_1 -  z_2||.$$
Hence, in view of (\ref{12})  the theorem is proved.

\begin{theorem} If $(t_0, c) \not \in S^+$ then the solution $z(t,t_0,c)$ of (\ref{6})   is unbounded on $[t_0,\infty).$
\label{t3}
\end{theorem} 
{\it Proof.} \rm
Assume, on the contrary,  that  $z(t) = z(t,t_0,z_0) = (u,v)$ is a bounded 
solution of (\ref{6}) and $(t_0,z_0) \not \in S^+.$ It is obvious that 

\begin{eqnarray}
&& u(t)  = U(t,t_0)u(t_0) + \int^t_{t_0}U(t,s)g_+(s,z(s),z(\beta(s)))ds,\nonumber\\
&& v(t) = V(t,t_0) \kappa  - \int^{\infty}_t V(t,s)g_- (s,z(s),z(\beta(s)))ds,
\label{13}
\end {eqnarray}
where
$$\kappa = v(t_0) + \int^{\infty}_{t_0} V(t,s)g_- (s,z(s),z(\beta(s)))ds,$$
and the improper integral converges and is bounded on $[t_0,\infty).$
But the inequality $||V(t,s)|| \geq K^{-1} \exp(\alpha(t-s)), t\geq s$ implies that $z(t)$ is bounded  only if
$\kappa =0.$ By Theorem \ref{t2}  $z(t)$ satisfies (\ref{6}) with $c=u(t_0).$ Hence $u(t_0),v(t_0),t_0$
satisfy (\ref{6}) and $(t_0,z_0) \not \in S^+.$ The contradiction proves our theorem.

It is not difficult to see that applying Lemma \ref{l2}   one can formulate and prove for the case  
$(- \infty, t_0]$ the theorems concerning the surface $S^-$ similar to the assertions for  $S^+.$

On the basis of  Theorems \ref{t1}-\ref{t3} and their analogues for 
$t \rightarrow - \infty,$ one can conclude that there exist two integral surfaces 
$\Sigma^+, \Sigma^-$ of equation (\ref{1}) such that every solution which starts at $\Sigma^+$ tends to  zero as 
$t\rightarrow \infty,$
and  every solution which starts at $\Sigma^-$ tends to  zero as $t\rightarrow -\infty.$
All solutions on $\Sigma^+, \Sigma^-$ are continued on $\mathbb R.$ If a solution  starts outside   $\Sigma^+$ then  it is unbounded on $[t_0, \infty),$  and if 
a solution  starts outside  $\Sigma^-$ then either it is unbounded on $(-\infty, t_0]$
or it cannot be back continued to $-\infty.$
\begin{remark}
Our results are not simply analogues of theorems for ordinary differential equations,  since we have proven the back continuation of solutions on the manifolds.  
\end{remark}
The following assertions of this section describe the structure of the set of solutions 
which do not belong to the integral sets.
Consider a solution $z(t) = \{u(t),v(t)\}$  of (\ref{1}) with the  initial data $t_0,u_0,v_0, z_0 \not = 0.$
\begin{lemma} 
Assume $(t_0,u_0,v_0) \not \in S^+,\,K(K^2 +1)(1+\exp(\alpha \theta))L < \sigma,$ and\\ $||u_0|| \leq ||v_0||.$  Then
$||u(t)|| \leq K^2||v(t)||, t \geq t_0.$
\label{l*}
\end{lemma}
{\it Proof.} Assume that  $K>1$ (otherwise the proof  is more simple).
Suppose on contrary that there exist  moments $\bar t, \tilde t$ such that 
\begin{eqnarray}
&& ||u(\bar t)||= ||v(\bar t)||,   ||u(\tilde t)||=K^2 ||v(\tilde t)||, 
\label{tt}
\end{eqnarray}
and  
\begin{eqnarray}
&&||v(t)|| \leq  ||u(t)|| \leq K^2 ||v(t)||,    \bar t <t  < \tilde t.
\label{ttt}
\end{eqnarray}
Since 
$$u(t) =   U(t,\bar t)u(\bar t) + \int^t_{\bar t}U(t,s)g_+(s,z(s),z(\beta(s)))ds,$$
we have that 
$$||u(t)|| \leq K e^{-\sigma(t- \bar t)}||u(\bar t)|| +  2KL(1+ e^{\alpha \theta}) \int^t_{\bar t} e^{-\sigma(t- s)}
    ||u(s)||ds.$$
Applying  Gronwell-Bellman Lemma  one can obtain that 
\begin{eqnarray}
&&||u(\tilde t)|| \leq K ||u(\bar t)||  e^{-(\sigma-   2KL(1+ e^{\alpha \theta})) (\tilde t- \bar t)}.
\label{ttt*}
\end{eqnarray}
Similarly we can find that
\begin{eqnarray}
&&||v(\bar t)|| \leq K ||v(\tilde t)||  e^{-(\sigma-   KL(1+K^2)(1+ e^{\alpha \theta})) (\tilde t- \bar t)}.
\label{ttt**}
\end{eqnarray}
Hence, using (\ref{ttt*}) and the equality $||u(\bar t)||= ||v(\bar t)||$ one can obtain
 the inequality $||u(\tilde t)||<K^2 ||v(\tilde t)||,$ which contradicts the assumption.
The Lemma is proved.
\begin{lemma} Assume that all conditions of Lemma \ref{l*}  are valid. Then 
\begin{eqnarray}
&&||v(t)|| \geq \frac{||v_0||}{K}e^{(\sigma-   KL(1+K^2)(1+ e^{\alpha \theta})) (t- t_0)}.
\label{ttt*****}
\end{eqnarray}
\label{ttt***}
\end{lemma}
{\it Proof.} Similarly to (\ref{ttt**}), applying Lemma \ref{l*} we may write that
$$||v(t_0)|| \leq K ||v(t)||e^{(\sigma-   KL(1+K^2)(1+ e^{\alpha \theta})) (t_0-t)}.$$
The last inequality is equivalent to (\ref{ttt*****}). The Lemma is proved.
\begin{theorem} 
Assume  that $(t_0,u_0,v_0) \not \in S^+$ and $K(K^2 +1)(1+ e^{\alpha \theta})L < \sigma.$ Then
$||v(t)|| \rightarrow \infty$ as $t \rightarrow \infty.$  
\label{t*}
\end{theorem}
{\it Proof.} If $||u_0|| \leq ||v_0||,$ then the proof is similar to that of Lemma \ref{ttt***}.
Assume that  $||u_0|| >||v_0||.$  In the same way as we obtained (\ref{ttt*}), we show that 
\begin{eqnarray}
&&||u(t)|| \leq K ||u(t_0)||  e^{-(\sigma-   2KL(1+ e^{\alpha \theta})) (t- t_0)}.
\label{ttt=}
\end{eqnarray}
Now, Theorem \ref{t3} and the last inequality imply  that there exists $ \bar t$ such that
$||u(\bar t)|| = ||v(\bar t)||.$   The theorem is proved.

Similarly to Theorem \ref{t*} one can prove  that the following theorem is valid.
\begin{theorem} 
Assume $(t_0,u_0,v_0) \not \in S^-$ and $K(K^2 +1)(1+ e^{\alpha \theta}))L < \sigma.$  Then
either  $||u(t)|| \rightarrow \infty$ as $t \rightarrow -\infty,$   or the solution 
$z(t)$ cannot be back continued to $-\infty.$ 
\label{t**}
\end{theorem}

\subsection{The smoothness of the surfaces}
The following condition is needed in this part of the paper.
\begin{itemize}
\item[C4)] The function $f(t,x,w)$ is uniformly continuously differentiable in $x,w$ for all
$t, x, w.$ 
\end{itemize}
\begin{theorem} Suppose conditions $C1)-C4)$ and inequality (\ref{bw}) are fulfilled. Then the function $F(t,u)$
is continuously differentiable in $u$ for all $t,u.$
\label{t4}
\end{theorem}
{\ Proof.}  By Theorem \ref{t1} there exists a solution  $z(t) = z(t,t_0,c)$ of (\ref{3*}).
Let us show that  it is continuously differentiable in $c.$
Denote \\$e_j = (0,\dots,0,1,0,\dots 0)^T, $ where the $j-$th coordinate is one, and  $h_j = he_j,$ where $h \in \mathbb R$ is  fixed.  
Denote
$\Delta z(t) =  z(t,t_0,c+h_j) - z(t,t_0,c), \Delta z = (\Delta u, \Delta v).$
We have that 
\begin{eqnarray}
&& \Delta u(t)  = U(t,t_0)h_j + \int^t_{t_0}U(t,s)\Big [\frac{\partial g_+(s,\Delta z(s),\Delta z(\beta(s)))}{\partial x} 
\Delta z(s) +\nonumber\\
&& \frac{\partial g_+(s,\Delta z(s),\Delta z(\beta(s)))}{\partial w} \Delta z(\beta(s)) +
o_+(\Delta z)\Big]ds,\nonumber\\
&&  \Delta v(t) = - \int^{\infty}_t V(t,s)\Big[\frac{\partial g_-(s, \Delta z(s),\Delta z(\beta(s)))}{\partial x}
 \Delta z(s)+\nonumber\\
&& \frac{\partial g_-(s,\Delta z(s),\Delta z(\beta(s)))}{\partial w} \Delta z(\beta(s)) +
o_-(\Delta z)\Big]ds.
\label{14}
\end {eqnarray}
Consider the following system of integral equations

\begin{eqnarray}
&& \phi(t)  = U(t,t_0)e_j + \int^t_{t_0}U(t,s)\Big [\frac{\partial g_+(s,\Delta z(s),\Delta z(\beta(s)))}{\partial x} 
\omega (s) +\nonumber\\
&& \frac{\partial g_+(s,\Delta z(s),\Delta z(\beta(s)))}{\partial w} \omega(\beta(s)) \Big]ds,\nonumber\\
&& v(t) =  - \int^{\infty}_t V(t,s) \Big[\frac{\partial g_-(s, \Delta z(s),\Delta z(\beta(s)))}{\partial x}
 \omega(s)+\nonumber\\
&& \frac{\partial g_-(s,\Delta z(s),\Delta z(\beta(s)))}{\partial w} \omega(\beta(s)) \Big]ds,
\label{15}
\end {eqnarray}
where  $\omega = (\phi, \psi).$

Condition $C4)$ implies that there exists a constant $L >0$ such that
$$ ||\frac{\partial g_+(s,\Delta z(s),\Delta z(\beta(s)))}{\partial x}|| \leq L,
 ||\frac{\partial g_+(s,\Delta z(s),\Delta z(\beta(s)))}{\partial w}|| \leq L,$$
$$ ||\frac{\partial g_-(s,\Delta z(s),\Delta z(\beta(s)))}{\partial x}|| \leq L,
 ||\frac{\partial g_-(s,\Delta z(s),\Delta z(\beta(s)))}{\partial w}|| \leq L.$$
Consequently, from Theorem \ref{t1} one can obtain  that there exists a unique bounded solution $\omega(t) =\omega(t,t_0,e_j) = (\phi,\psi)$
of system (\ref{15}).  Equalities  (\ref{14}) and (\ref{15}) imply that 
\begin{eqnarray}
&& \frac{ \Delta u(t)}{h}-\phi  = 
\int^t_{t_0}U(t,s)\Big [\frac{\partial 
g_+(s,\Delta z(s),\Delta z(\beta(s)))}{\partial x}\Big( \frac{ \Delta z(s)}{h}-\omega(s)\Big) +\nonumber\\
&& \frac{\partial g_+(s,\Delta z(s),\Delta z(\beta(s)))}{\partial w}\Big(\frac{ \Delta z(\beta(s))}{h}-\omega(\beta(s)\Big) +   
\frac{o_+(\Delta z)}{h}\Big]ds,\nonumber\\
&&   \frac{ \Delta v(t)}{h}-\psi   = - \int^{\infty}_t V(t,s)\Big[\frac{\partial g_-(s, \Delta z(s),\Delta z(\beta(s)))}{\partial x}
 \Big( \frac{ \Delta z(s)}{h}-\omega(s)\Big)    +\nonumber\\
&& \frac{\partial g_-(s,\Delta z(s),\Delta z(\beta(s)))}{\partial w}\Big( \frac{ \Delta z(s)}{h}-\omega(s)\Big) +
\frac{o_-(\Delta z)}{h}\Big]ds.
\label{16}
\end {eqnarray}
Since solutions of (\ref{6}) satisfy Lipshitz condition, $||\Delta z(t)|| \leq 2 K |h| \exp(-\alpha(t-t_0))$ and  $C4)$ is valid, 
$$||o_\pm|| \leq O(h) |h|  \exp(-\alpha(t-t_0)), i = 1,2.$$
Using the last relation and (\ref{16}) we can check that 
$$\sup_{[t_0,\infty)} || \frac{ \Delta z}{h}-\omega|| \leq \frac{2KL\exp(\sigma)}{ \alpha}
\sup_{[t_0,\infty)} || \frac{ \Delta z}{h}-\omega|| +  \frac{2KL\exp(\sigma)}{ \alpha} O(h),$$
where $O(h) \rightarrow 0$  as $h \rightarrow 0.$
By the assumption we can take $L$ arbitrarily  small. Then  $$\frac{\partial z(t,t_0,c)}{\partial c_j} = \omega(t).$$ The theorem is proved.

The last theorem, and its analogue for $S^-,$ which can be proved similarly, imply the smoothness of the surfaces $\Sigma^+$ and $\Sigma^-.$
\subsection{The existence of  bounded and  periodic solutions}
In the same way as we did in  Lemma \ref{l1}, we can prove that the following assertion is valid.

\begin{lemma} Assume that conditions $C1)-C3)$
and inequality (\ref{bw}) are valid.  
A bounded on $\mathbb R$ function $z(t) = (u,v)$ is a solution of (\ref{3})  if and only if it is a solution of the following
system of integral equations
\begin{eqnarray}
&& u(t) =  \int^t_{-\infty}U(t,s)g_+(s,z(s),z(\beta(s)))ds,\nonumber\\
&& v(t) = - \int^{\infty}_t V(t,s)g_- (s,z(s),z(\beta(s)))ds.
\label{6*}
\end{eqnarray}
\label{ll}
\end{lemma} 
In what follows, we shall need the  following  two conditions:
\begin{itemize}
\item[C5)]  $||f(t,0,0)|| < h, t \in \mathbb R, {\mbox for \, some} \,  h \in \mathbb R, h>0;$
\item[C6)]  $\frac{2KL}{\sigma}<1.$
\end{itemize}
Denote  $H= \sup_{\mathbb R} ||U(t)|| h.$ 
\begin{theorem} Suppose conditions $C1)-C3), C5),C6)$ and inequality (\ref{bw}) are valid.
Then there 
exists a unique bounded on $\mathbb R$ solution of  (\ref{3}), $z(t),$  and 
\begin{eqnarray}
&& ||z(t)|| \leq \frac{2KH}{\sigma -2KL}.
\label{7****}
\end{eqnarray}
\label{t1*}
\end{theorem}
{\it Proof.} Let us solve  system (\ref{6*}), applying  the method of successive approximations.
Assume that $z_0(t) = 0, t \in \mathbb R,$  and   
\begin{eqnarray}
&& u_m(t) =  \int^t_{-\infty}U(t,s)g_+(s,z_{m-1}(s),z_{m-1}(\beta(s)))ds,\nonumber\\
&& v_m(t) = - \int^{\infty}_t V(t,s)g_- (s,z_{m-1}(s),z_{m-1}(\beta(s)))ds,
\label{6**}
\end{eqnarray}
where $z_m(t) = (u_m,v_m)^T, m \geq 1.$
We have that
$$||u_1(t)|| \leq  \int^t_{-\infty}||U(t,s)|| ||g_+(s,0,0)||ds, $$$$||v_1(t)|| \leq  \int^{\infty}_t ||V(t,s)|| ||g_- (s,0,0)||ds.$$
Then 
$$ ||u_1(t)|| \leq \frac{KH}{\sigma},  $$$$||v_1(t)|| \leq \frac{KH}{\sigma},$$
and hence  $$||z_1(t)|| \leq   \frac{2KH}{\sigma}.$$
We can obtain by induction  that the following
inequality is valid
$$||z_{m+1}(t) - z_{m}(t)|| \leq 
\Big( \frac{2KL}{\sigma}\Big)^{m+1}\frac{H}{L}.$$ Condition $C6)$
implies that the sequence
$z_m$  converges uniformly to a function $\kappa(t), t \in \mathbb R.$ It is easy to check that
$$||\kappa(t) || \leq \frac{2KH}{\sigma -2KL}$$ and $\kappa(t)$ is a solution
of system (\ref{6*}). Assume  that there are two different bounded solutions $z^1, z^2$ of (\ref{6*}).
Then one can see that
$$\sup_{\mathbb R} ||z^1(t) - z^2(t)|| \leq  \frac{2KL}{\sigma}\sup_{\mathbb R} ||z^1(t) - z^2(t)||.$$
The last formula contradicts  assumption $C6).$ The theorem is proved.

Let us make the following additional assumptions:
\begin{itemize}
\item[C7)] the matrix $A(t)$ and the function $f(t,x,y)$ are periodic in $t$ with a period  $\omega;$
\item[C8)]there exist  a real number $\bar \omega$ and a positive integer $p,$ such that
$\theta_ {i+p} =\theta_ {i}+ \bar \omega, i \in \mathbb Z.$ 
\end{itemize}
\begin{theorem} Suppose conditions $C1)-C3),$ and $ C5)-C8)$ are valid. Moreover,
$$\frac{\omega}{\bar \omega} = \frac{k}{m}, k,m \in \mathbb N.$$  
Then there 
exists an $m \omega-$ periodic solution of (\ref{3}).
\end{theorem}
Indeed, one can verify that all approximations $z_m(t)$ are $m \omega-$ periodic
functions. The definition of $\kappa(t)$ implies that it is $m \omega-$ periodic, too.
The theorem is proved.

\begin{remark}
The last theorem can be applied to equations of logistic type \cite{s2}.
Another  interesting possibility is to consider the existence of quasiperiodic solutions using the method
of equivalent integral equations \cite{y1} . 
\end{remark}

\noindent {\bf e-mail:} marat@metu.edu.tr
\end{document}